\documentclass[preprint,12pt,authoryear]{elsarticle}

\usepackage{mathrsfs}
\usepackage{xcolor}
\usepackage{mathscinet}
\usepackage{latexsym}
\usepackage{amsthm}
\usepackage{amssymb}
\usepackage{amsfonts}
\usepackage{amsmath}
\usepackage{natbib}
\usepackage{graphicx}
\usepackage{float}

\newtheorem{theorem}{Theorem}[section]
\newtheorem{thm}[theorem]{Theorem}
\newtheorem{lemma}[theorem]{Lemma}
\newtheorem{lem}[theorem]{Lemma}
\newtheorem{proposition}[theorem]{Proposition}

\newtheorem{corollary}[theorem]{Corollary}

\theoremstyle{definition}

\newtheorem{defn}[theorem]{Definition}

\theoremstyle{remark}
\newtheorem{remark}[theorem]{Remark}
\newtheorem{rem}[theorem]{Remark}
\numberwithin{equation}{section}

\newcommand{\E}{\mathbb{E}}            

\newcommand{\R}{\mathbb{R}}

\newcommand{\PP}{\mathbb{P}}

\newcommand{\Be}{\begin{equation}}
\newcommand{\Ee}{\end{equation}}
\newcommand{\Bs}{\begin{split}}
	\newcommand{\Es}{\end{split}}
\newcommand{\Bes}{\begin{equation*}}
	\newcommand{\Ees}{\end{equation*}}
\newcommand{\BT}{\begin{thm}}
	\newcommand{\ET}{\end{thm}}
\newcommand{\Bp}{\begin{proof}}
	\newcommand{\Ep}{\end{proof}}
\newcommand{\BL}{\begin{lem}}
	\newcommand{\EL}{\end{lem}}
\newcommand{\BP}{\begin{proposition}}
	\newcommand{\EP}{\end{proposition}}
\newcommand{\BC}{\begin{corollary}}
	\newcommand{\EC}{\end{corollary}}
\newcommand{\BR}{\begin{rem}}
	\newcommand{\ER}{\end{rem}}
\newcommand{\BD}{\begin{defn}}
	\newcommand{\ED}{\end{defn}}
\newcommand{\BI}{\begin{itemize}}
	\newcommand{\EI}{\end{itemize}}

\newcommand{\varX}{\boldsymbol{X}}
\newcommand{\designmatrix}{\boldsymbol{X}}
\newcommand{\designrow}{\boldsymbol{x}}
\newcommand{\designrowi}{\boldsymbol{x}_i}
\newcommand{\designs}{x}
\newcommand{\designsij}{x_{ij}}
\newcommand{\varY}{Y}
\newcommand{\outcome}{\boldsymbol{Y}}
\newcommand{\outcomes}{y}
\newcommand{\outcomesi}{y_i}

\newcommand{\parameter}{\boldsymbol{\beta}}
\newcommand{\target}{\boldsymbol{\beta}^*}
\newcommand{\est}{\widehat{\boldsymbol{\beta}}}
\newcommand{\penaltyest}{\widehat{\lambda}}

\newcommand{\varW}{\boldsymbol{W}} 
\newcommand{\varWi}{\varW_i}
\newcommand{\varWs}{W}
\newcommand{\varWij}{W_{ij}}
\newcommand{\sumW}{\boldsymbol{S}^{\varWs}_{n}}
\newcommand{\sumWj}{S^{\varWs}_{n,j}}
\newcommand{\varWstilde}{\tilde{W}}
\newcommand{\varWitilde}{\tilde{\varW}_i}
\newcommand{\sumWtilde}{\boldsymbol{S}^{\varWstilde}_{n}}

\newcommand{\sumWltilde}{S^{\varWstilde}_{n,l}}

\newcommand{\varZ}{\boldsymbol{Z}} 
\newcommand{\varZi}{\varZ_i}
\newcommand{\varZs}{Z}
\newcommand{\varZij}{Z_{ij}}
\newcommand{\sumZ}{\boldsymbol{S}^{\varZs}_{n}}
\newcommand{\sumZj}{S^{\varZs}_{n,j}}
\newcommand{\varZstilde}{\tilde{Z}}
\newcommand{\varZitilde}{\tilde{\varZ}_i}
\newcommand{\sumZtilde}{\boldsymbol{S}^{\varZstilde}_{n}}

\newcommand{\sumZltilde}{S^{\varZstilde}_{n,l}}

\newcommand{\fcnf}{f}

\journal{}

\begin{document}

\begin{frontmatter}



\title{Estimating the Penalty Level of $\ell_1$-minimization via Two Gaussian Approximation Methods}


\author{Fang Xie}
\address{School of Mathematics and Statistics, Wuhan University,\\
	Wuhan, Hubei 430072, P.R. China.\\}
\ead{fangxie219@foxmail.com}

\begin{abstract}
In this paper, we aim to give a theoretical approximation for the penalty level of $\ell_{1}$-regularization problems. 
This can save much time in practice compared with the traditional methods, such as cross-validation. 
To achieve this goal, we develop two Gaussian approximation methods, which are based on a moderate deviation theorem and Stein's method respectively. 
Both of them give efficient approximations and have good performances in simulations. 
We apply the two Gaussian approximation methods into three types of ultra-high dimensional $\ell_{1}$ penalized regressions: lasso, square-root lasso, and weighted $\ell_1$ penalized Poisson regression. 
The numerical results indicate that our two ways to estimate the penalty levels achieve high computational efficiency. 
Besides, our prediction errors outperform that based on the 10-fold cross-validation. 
\end{abstract}



\begin{keyword}


penalty level\sep Gaussian approximation\sep Stein's method\sep moderate deviation theorem\sep generalized linear model
\end{keyword}

\end{frontmatter}


\section{Introduction}\label{sec:Introduction}
With the increase of the dimension of data, high-dimensional regressions are widely applied in many areas, such as economics~\citep{Bai2008Forecasting}, biology~\citep{Xie2019Aggregating}, health science~\citep{Riphahn2003Incentive}. 
At the same time, it becomes more and more challenging when the dimension of variables, $p$, is much larger than the sample size, $n$. 
Either for linear regression or generalized linear regression, variable selection is a good way to reduce the dimension when $p$ is large. 
A family of the popular variable selection methods is based on the $\ell_1$ regularization. 
For instance, lasso \citep{Tibshirani1996Regression}, square-root lasso \citep{Belloni2011Square-root}, and $\ell_1$ penalized Poisson regression \citep{Li2015Consistency}. 
In these methods, choosing a suitable penalty level is very important since it can influence the estimation accuracy directly. 
In general, the penalty level is selected by cross-validation, C$_p$, AIC or BIC criterions \citep{Chen2008Extended,Tibshirani1996Regression, Zou2007On}. 
They can produce good estimation accuracy, but heavy procedures are involved. 
In this paper, we give two theoretical approximations for the penalty level, and they can reduce much computation cost in practice and improve the prediction accuracy at the same time.

Consider the following generalized linear models \citep{mccullagh1989generalized} which include the classical linear models. 
Suppose $(\designrow_1,\outcomes_1),\ldots,(\designrow_n,\outcomes_n)$ are independent pairs of observed data which are realizations of random vectors $(\varX_{1},\varY_{1}),\ldots,(\varX_{n},\varY_{n})$,
with $p$-dimensional covariates $\varX_{i}\in \R^p$ and univariate response variables $\varY_{i}\in \R$ for all $i\in\{1,\ldots,n\}$. 
$(\varX_{i},\varY_{i})$ are assumed to satisfy the conditional distribution
\begin{equation}\label{e:glm1}
\varY_{i}|\varX_{i}=\designrowi\sim F{\rm\ \ \ \ with}\ \ \ \ g(\E(\varY_{i}|\varX_{i}=\designrowi))=\designrowi'\target,
\end{equation}
where $F$ is a distribution in the exponential family, $g(\cdot)$ a real-valued link function, $\E(\cdot)$ the expectation function, and $\target\in\mathbb{R}^{p}$ an unknown parameter vector. 
We denote by $A'$ the transpose of $A$ (a vector or matrix). 
If $F$ is the normal distribution and $g(t)=t$, then \eqref{e:glm1} is the classical linear model. 
If $F$ is the Poisson distribution and $g(t) = \log{t}$ for all $t>0$, then~\eqref{e:glm1} is the Poisson regression model. 

Denoting $\designrowi=(\designs_{i1},\cdots,\designs_{ip})'$, without loss of generality, we assume
$$\frac1{n}\sum\limits_{i=1}^n{\designsij}=0\ \ \ \ \ \ \ \ {\rm and}\ \ \ \ \ \ \ \frac1{n}\sum\limits_{i=1}^n{\designsij^2 }=1,\ \ \ \ \ \ \ \ \ \ \ {\rm for\ all}\ j\in[p].$$
Write $\designmatrix=(\designrow_1,\ldots,\designrow_n)'\in\R^{n\times p}$ and $\outcome=(\outcomes_1,\ldots,\outcomes_n)'\in \R^n$. 
Consider the following $\ell_{1}$-minimization problem:
\Be\label{e:betahat}
\est = \arg\min_{\parameter\in\R^p}\left\{ L(\parameter|\designmatrix,\outcome)+\lambda \|\parameter\|_1 \right\},
\Ee
where $L(\parameter|\designmatrix,\outcome)$, connected to the distribution of $\outcome|\designmatrix$, is assumed to be a convex function with respect to $\parameter$, and $\lambda>0$ the penalty level to be chosen. 
For example, 
if $F$ is normal distribution and $L(\parameter|\designmatrix,\outcome) = \|\outcome-\designmatrix\parameter\|_2^2/2n$, $\est$ in~\eqref{e:betahat} is the lasso estimator~\citep{Tibshirani1996Regression}; 
if $F$ is normal distribution and $L(\parameter|\designmatrix,\outcome) =  \sqrt{\|\outcome-\designmatrix\parameter\|_2^2/n}$, $\est$ in~\eqref{e:betahat} is the square-root lasso estimator~\citep{Belloni2011Square-root}. 

As we mentioned above, the choice of $\lambda$ affects the accuracy of $\est$ directly. 
In the previous research, the $\lambda$ had been proved to be well approximated by a factor times a Gaussian quantile in the case of the linear model with independent Gaussian errors. 
The lasso estimator can achieve the near-oracle performance with probability approaching to $1-\alpha$, if $\lambda=c\sigma(\sqrt{n})^{-1}\Phi^{-1}(1-{\alpha/2p})$ with $c>1$, $\sigma$ being the standard deviation of Gaussian error, $\Phi^{-1}(\cdot)$ the inverse of the cumulative distribution function of standard normal distribution, and $\alpha\in(0,1)$~\citep{Bickel2009Simultaneous}.
The square-root lasso was proved to own the same property~\citep{Belloni2011Square-root}, if $\lambda = c(\sqrt{n})^{-1}\Phi^{-1}(1-\alpha/2p)$ with $c>1$. 
But for the generalized linear model except for the linear model, there are a few papers to approximate the penalty theoretically. 
So, we develop two Gaussian approximation methods to estimate the penalty level for the  generalized $\ell_{1}$-minimization problem. 

We briefly introduce how to connect the two Gaussian approximation methods to the estimation of $\lambda$ and the corresponding results. 
Using Karush-Kuhn-Tucker conditions for \eqref{e:betahat}, we have
\Bes
\nabla L(\est|\designmatrix,\outcome) + \lambda \boldsymbol{\kappa} = 0,
\Ees
where $\nabla$ is the gradient operator, and $\kappa$ is the subgradient of $\|\parameter\|_1$ at $\parameter = \est$. 
Noticing $\|\boldsymbol{\kappa}\|_\infty\le 1$, the equality above implies 
\Bes
\lambda\ge\|\nabla L(\est|\designmatrix,\outcome)\|_\infty,
\Ees
where $\|\cdot\|_{\infty}$ is the $\ell_{\infty}$ norm.
Naturally, an ideal choice of penalty level $\lambda$ should guarantee that
\begin{equation*}
	\lambda\ge c\|\nabla L(\target|\designmatrix,\outcome)\|_{\infty},
\end{equation*}
where $c>1$ is a tuning parameter (see \cite{Belloni2011Square-root,Jia2019Sparse}).
Thus given a small $\alpha \in (0,1)$, we need to find a suitable $\lambda$ such that
\Be  \label{e:Hquantile}
\PP\left(c\|\nabla L(\target|\designmatrix,\outcome)\|_{\infty} \le \lambda\right)\ge 1-\alpha.
\Ee
We observe that for the generalized linear models, $\nabla L(\target|\designmatrix,\outcome)$ usually has a special form $\nabla L(\target|\designmatrix,\outcome)=\sum_{i=1}^n \fcnf(\designrowi,\outcomesi)/n$ with $\fcnf(\designrowi,\outcomesi)$ being independent random variables depending on the specific models (see more details in Sections~\ref{sec:GauAppr} and~\ref{sec:Applications}). 
For simplicity, we denote $\varWi=\fcnf(\designrowi,\outcomesi)$.

In this paper, we propose two approximated penalty levels for $\lambda$. 
One is $\penaltyest_1=c\theta(\sqrt{n})^{-1}\Phi^{-1}(1-\alpha/2p)$ with $c>1$, $\theta\in(0,+\infty)$ and $\alpha\in(0,1)$, under which we prove that 
\begin{equation*}
	\PP\left(c\|\nabla L(\target|\designmatrix,\outcome)\|_{\infty} \le \penaltyest_1\right)\ge 1-\alpha(1+O(\frac{(\log{p})^{5/2}}{\sqrt{n}})).
\end{equation*}
The result is obtained by utilizing a moderate deviation theorem that we state in Section~\ref{sec:GauAppr}. 
This Gaussian approximation was also used in the previous research, see \citep{Bickel2009Simultaneous} for lasso, \citep{Belloni2011Square-root} for square-root lasso. 
But they assumed that $F$ was a normal distribution. 
In our case, we allow $F$ to be any distribution in the exponential family. 
The other is $\penaltyest_2=c(\sqrt{n})^{-1}z_{1-\alpha}$ with $z_{1-\alpha}$  satisfying 
\begin{equation*}
	\PP\left(\max_{j\in[p]\}}\left|\frac{1}{\sqrt{n}}\sum_{i=1}^{n}\varZij\right|\le z_{1-\alpha}\right)=1-\alpha,
\end{equation*}
and $\varZi=(\varZs_{i1},\ldots,\varZs_{ip})'\sim N(0,\E[\varWi\varWi'])$. 
Under this approximated penalty level, we prove by Stein's method that 
\begin{equation*}
	\PP\left(c\|\nabla L(\target|\designmatrix,\outcome)\|_{\infty} \ge \penaltyest_2\right)\le 1-\alpha-O(n^{-1/8}(\log{p})^{7/8}).
\end{equation*}

Under both two approximations, we see that the inequality~\eqref{e:Hquantile} holds when $n,p\rightarrow\infty$. 
The difference is in the first case, $n,p$ need to satisfy $p\le e^{o(n^{1/5})}$, and a stronger condition $p\le e^{o(n^{1/7})}$ has to be assumed in the second case. 
The more detailed comparisons for the two cases are referred to Section~\ref{sec:GauAppr}.

The rest of the paper is organized as follows. 
Section \ref{subsec:notations} below gives the notations throughout the paper. 
In Section \ref{sec:GauAppr}, we give two theoretical approximations of penalty level $\lambda$ and the corresponding proofs are deferred to Appendix. 
In Section \ref{sec:Applications}, we apply the two approximations to three types of $\ell_1$ penalized regression. 
Besides, we conduct simulations to show the prediction errors and computation time of the three methods under the two approximated penalty levels and compare them with the results of 10-fold cross-validation. 
Finally, we give a conclusion in Section~\ref{sec:conclusion}. 

\subsection{Notations}\label{subsec:notations}
For simplicity of notations, we use $\E_{n}(\cdot)$ denotes the average over index $i\in [n]$, where $[n] = \{1,2,\ldots,n\}$. 
For example, $\E_n(\cdot) = (\sum_{i=1}^n (\cdot))/n$. 
$a\lesssim b$ means that it exists a universal positive constant $c$ such that $a\le c b$.
For a $d$-dimensional vector $\boldsymbol{v}=(v_1,\ldots,v_d)'$, we denote by $\|\boldsymbol{v}\|_q$ its $l_q$ norm for all $q\ge 1$. 
Especially, when $q=\infty$, $\|\boldsymbol{v}\|_{\infty} = \max\limits_{i\in[d]}|v_i|$. 
The notation $b_n=O(a_n)$ implies $b_n\lesssim a_n$, and $b_n = o(a_n)$ implies $\lim\limits_{n\rightarrow\infty}b_n/a_n=0$. 
Especially, $O(1)$ stands for a positive finite constant and $o(1)$ is infinitesimal.

\section{Two Gaussian approximation methods to estimate the penalty level}\label{sec:GauAppr}
In this section, we introduce the two Gaussian approximation methods for estimating the penalty level of $\ell_{1}$-regularized regression. 
One is based on a moderate deviation theorem, the other is based on Stein's method. 

We recall that the good approximation of penalty level $\lambda$ should satisfy for given $\alpha\in(0,1)$, 
\begin{equation}\label{ineq:ideallambda}
\PP\left(c\|\nabla L(\target|\designmatrix,\outcome)\|_{\infty} \le \lambda\right)\ge 1-\alpha,
\end{equation}
with a suitable constant $c>1$. 
Observe that for the generalized linear model, $\nabla L(\target|\designmatrix,\outcome)$ can be written as 
\begin{equation*}
	\nabla L(\target|\designmatrix,\outcome)=\frac1n\sum_{i=1}^n\varWi,
\end{equation*}
with $\varWi=\fcnf(\designrowi,\outcomesi)\in\R^p$ being independent random vectors. 
For example, $\varWi=\designrowi(\designrowi'\target-\outcomesi)$ for the lasso; $\varWi=\designrowi(\designrowi'\target-\outcomesi)/\sqrt{\|\outcome-\designmatrix\target\|_2^2/n}$ for the square-root lasso. 

Denote $\varWi=(\varWs_{i1},\ldots,\varWs_{ip})'$ and 
\begin{equation*}
	\sumWj = \frac1{\sqrt{n}}\sum_{i=1}^n \varWij.
\end{equation*}
Then, we aim to prove that for a suitable $\lambda$, the following probability
\begin{equation*}
	\PP\left(c\|\nabla L(\target|\designmatrix,\outcome)\|_{\infty}\le\lambda\right)=\PP\left(\max_{j\in [p]} |\sumWj|\le \lambda\sqrt{n}/c\right)
\end{equation*}
is close to 1. 

\subsection{Gaussian approximation based on moderate deviation theorem}
In this section, we prove that $\penaltyest_1=c\theta(\sqrt{n})^{-1}\Phi^{-1}(1-\alpha/2p)$ with $c>1$, $\theta\in(0,+\infty)$ and $\alpha\in(0,1)$ is a good approximation of $\lambda$ so that the inequality~\eqref{ineq:ideallambda} holds. 
The main skill is the moderate deviation theorem that was given by \citep{Sakhanenko1991}. 
One version stated in Lemma~\ref{lem:MDT} is very convenient for use. 
\begin{lem} [\citep{Liu2013Self-normalized}]
	\label{lem:MDT}  Let $\eta_{1},\cdots,\eta_{n}$ be independent random
	variables with $\mathbb{E}\eta_{i}=0$ and $|\eta_{i}|\le 1$ for all $i\in[n]$.
	Denote $\sigma_{n}^{2}=\sum\limits_{i=1}^{n}\mathbb{E}\eta_{i}^{2}$ and $%
	T_{n}=\sum\limits_{i=1}^{n}\mathbb{E}|\eta_{i}|^{3}/\sigma_{n}^{3}$. Then
	there exists a positive constant $K$ such that for all $x\in[%
	1,\min\{\sigma_{n},L_{n}^{-1/3}\}/K]$
	\begin{equation*}
		\mathbb{P}(\sum\limits_{i=1}^{n}\eta_{i}>x\sigma_{n})=(1+O(1)x^{3}T_{n})\bar%
		\Phi(x),
	\end{equation*}
	where $\bar\Phi(x)=1-\Phi(x)$ and $\Phi(x)$ is the cumulative distribution
	function of standard normal distribution.
\end{lem}
This type of moderate deviation theorem was also applied in \citep{Jia2019Sparse,Liu2013Self-normalized}. 
Other type moderate deviation theorems can refer to \citep{Hu2009Cramer,Jing2008Towards,Liu2010Cramer,Shao2016Cramer}.

The following theorem is the key to proving that $\penaltyest_1$ is a good approximation of $\lambda$. 
It gives the upper bound of $\PP\left(\max_{j\in[p]} |\sumWj| \le z\right)$ and is proved by means of a truncation technique and Lemma \ref{lem:MDT}, see~\ref{appendix:MDT} below.
\begin{theorem}\label{thm:TailProbMDT}
	Suppose for each $i\in[n],j\in[p]$, $\E \varWij=0$ and $\E_n(\E \varWij^2)=\theta^2<\infty$. 
	Assume that $\sup_{i\in[n],j\in[p]}\E e^{t_1|\varWij|}<\infty$ for some $0<t_1<\infty$. 
	Then for all $z\in [C\sqrt{\log{p}},o(n^{1/6}(\log{p})^{-1/3}))$ with a positive constant $C$, we have
	\begin{equation}\label{ineq:TailProbMDT}
	\begin{split}
	&\PP\left(\max_{j\in[p]} |\sumWj| \le z\right)\\
	&\ge 1-2p\bar{\Phi}\left(z\right)\left(1+O(1)\frac{(z-(\log{p})/p^3)^3\log{p} }{\sqrt{n}}\right)\left(1+\frac{O(1)\log{p}}{p^3\bar{\Phi}\left(z\right)}\right) - Cn/p^2,
	\end{split}
	\end{equation}
	where $\bar{\Phi}(x) = 1- \Phi(x)$ and $\Phi(x)$ is the cumulative distribution function of standard normal distribution.
\end{theorem}

\begin{remark}\label{rem:MDT}
	Noticing the range of $z$, \eqref{ineq:TailProbMDT} can be simplified to
	$$\PP\left(\max_{j\in[p]} |\sumWj| \le z\right)\ge 1-2p\bar{\Phi}\left(z\right)(1+O(\frac{(\log p)^{5/2}}{\sqrt{n}})),$$
	when $p$ is sufficiently large.
\end{remark}

By the fact of $\penaltyest_1\sqrt{n}/c\sim\sqrt{\log{p}}$ and Theorem~\ref{thm:TailProbMDT}, we can obtain the following corollary directly. 
\begin{corollary}\label{coro:MDT}
	Suppose that $\nabla L(\target|\designmatrix,\outcome)$ can be written as the following form
	\begin{equation*}
		\nabla L(\target|\designmatrix,\outcome)=\frac1n\sum_{i=1}^n\varWi,
	\end{equation*}
	with $\varWi=\fcnf(\designrowi,\outcomesi)\in\R^p$ being independent random vectors. 
	Assume that for each $i\in[n],j\in[p]$, $\E \varWij=0$, $\E_n(\E \varWij^2)= \theta^2<\infty$ and  $\sup_{i\in[n],j\in[p]}\E e^{t_1|\varWij|}<\infty$ for some $0<t_1<\infty$. 
	Then, we have 
	\begin{equation*}
		\PP\left(c\|\nabla L(\target|\designmatrix,\outcome)\|_{\infty}\le\penaltyest_1\right)\ge 1-\alpha(1+O(\frac{(\log p)^{5/2}}{\sqrt{n}})),
	\end{equation*}
	with $\penaltyest_1=c\theta(\sqrt{n})^{-1}\Phi^{-1}(1-\alpha/2p)$, $c>1$ and $\alpha\in(0,1)$. 
	Furthermore, when $n, p\rightarrow\infty$ obeying $p\le e^{o(n^{1/5})}$, we have 
	\begin{equation}\label{ineq:lambda1}
	\PP\left(c\|\nabla L(\target|\designmatrix,\outcome)\|_{\infty}\le\penaltyest_1\right)\ge 1-\alpha(1+o(1)),
	\end{equation}
\end{corollary}

\begin{rem}\label{rem:coroMDT}
	Let the function $L$ is the negative log-likelihood of generalized linear model, and suppose $\outcomesi=\mu(\designrowi'\target) + \epsilon_i$ for all $i\in[n]$ with $\mu(\designrowi'\target)=g^{-1}(\designrowi'\target)$. 
	Then, the gradient of $L$ has the form 
	\begin{equation*}
		\nabla L(\target|\designmatrix,\outcome) = \frac1n\sum_{i=1}^n \designrowi\epsilon_i,
	\end{equation*}
	So, $\varWi$ in Corollary~\ref{coro:MDT} is equal to $\designrowi\epsilon_i$. 
	Only if we assume $\epsilon_i$ satisfies $\E\epsilon_i=0$,  $\E\epsilon_i^2=\sigma^2<\infty$, and $\sup_{i\in[n],j\in[p]}\E e^{t_1|\designsij\epsilon_i|}<\infty$ for some $0<t_1<\infty$, all the results of Corollary~\ref{coro:MDT} hold. 
\end{rem}

\subsection{Gaussian approximation based on Stein's method}
In this section, we aim to prove that $\penaltyest_2=c(\sqrt{n})^{-1}z_{1-\alpha}$ with $c>1$, $\alpha\in(0,1)$ and $z_{1-\alpha}$ satisfying 
\begin{equation*}
	\PP\left(\max_{j\in[p]\}}\left|\frac{1}{\sqrt{n}}\sum_{i=1}^{n}\varZij\right|\le z_{1-\alpha}\right)=1-\alpha,
\end{equation*}
is another good approximation of $\lambda$ so that the inequality~\eqref{ineq:ideallambda} holds. 
The main skills include Stein's method and truncation technique. 

Firstly, we give some assumptions and notations. 
Let $\varW_1,...,\varW_n$ be a sequence of independent $p$-dimensional random vectors with the following assumptions: $E(\varWi)=\boldsymbol{0}_p$ for $i=1,...,n$ and
$$Q_{jk}=\frac 1n \sum_{i=1}^n \E[\varWij \varWs_{ik}], \ \ \ \ \ j,k=1,...,p.$$
Obviously, $\boldsymbol{Q}=(Q_{jk})_{1 \le j,k \le p}$ is a $p \times p$ symmetric matrix. 
Define
$$\sumW=\frac{1}{\sqrt n}\sum_{i=1}^n \varWi\ \ \ {\rm and}\ \ \ \sumWj = \frac1{\sqrt{n}}\sum_{i=1}^n \varWij.$$
Let $\varZ_1,...,\varZ_n$ be a sequence of independent $p$-dimensional Gaussian random vectors such that $\varZi \sim N(\boldsymbol{0}_p,\E[\varWi\varWi'])$
for $i=1,...,n$. 
We assume $\{\varWi\}_{i=1}^n$ and $\{\varZi\}_{i=1}^n$ are independent. Denote
$$\sumZ=\frac{1}{\sqrt n}\sum_{i=1}^n \varZi\ \ \ {\rm and}\ \ \ \sumZj= \frac1{\sqrt{n}}\sum_{i=1}^n \varZij$$
and then we have $\sumZ \sim N(\boldsymbol{0}_p,\boldsymbol{Q})$.

\cite{chernozhukov2013Gaussian} put forward an ultrahigh dimensional Gaussian approximation by Stein's method, which gave a Berry-Esseen bound between $\max_{j\in[p]}\sumWj$ and $\max_{j\in[p]}\sumZj$.
The results have been applied to bootstrap and Dantzig. 
In this paper, we extend its application to estimate the penalty level $\lambda$ of $\ell_{1}$-regularized regressions. 
For more details about Gaussian approximation by Stein's method, we refer the reader to \cite{Chatterjee2005AnEB,Chen2011Nornal,Chen2011Multi,Chen2004Normal,chernozhukov2014Gaussian,CHERNOZHUKOV2016Empirical,Chernozhukov2017Central,Rollin2013Stein,stein1981Estimation}.

Follow the notations in \citep{chernozhukov2013Gaussian} and recall $\E_{n}(\cdot)$ denotes the average over index $1 \le i \le n$, that is, it simply abbreviates the notation $(\sum_{i=1}^{n}(\cdot))/n$. 
For instance, $\E_{n}(x_{ij}^{2})=\sum_{i=1}^{n} x^{2}_{ij}/n$. 
We define
\begin{equation*}
	a_\varWs(\gamma) = \inf\left\{a\ge 0: \PP\left(|\varWij|\le a\sqrt{\E_n(\E\varWij^2)}, {\rm \ for\ all\ } i\in [n], j\in [p]\right)\ge 1-\gamma\right\},
\end{equation*}
where $\gamma\in (0,1)$.
Similarly, we can define $a_\varZs(\gamma)$ and then define
$$a(\gamma) = a_\varWs(\gamma)\vee a_\varZs(\gamma),$$
where the notation $c\vee d = \max\{c, d\}$. 
Denote
\begin{equation*}
	M_k = \max_{j\in [p]}\left(\E_{n}\left(\E\left|\varWij\right|^k\right)\right)^{1/k}.
\end{equation*}

The following theorem can be derived by Theorem 2 and Lemma 2 of \cite{chernozhukov2013Gaussian} and some simple transformations. 
The proof is given in~\ref{appendix:SM}.
\begin{theorem} \label{thm:GauAppr}
	Suppose $c_1\le M_2^2\le C_1$ and $\E_n(\E\max_{j\in [p]} \varWij^4)\le C_2$ with some positive constants $c_1$, $C_1$ and $C_2$. 
	For any $\gamma\in(0,1)$, we have
	\begin{equation} \label{e:ULB}
	\begin{split}
	&\sup_{z\ge 0}\left|\PP\left(\max_{j\in[p]} |\sumWj| \le z\right)-\PP\left(\max_{j\in[p]} |S^Z_{n,j}| \le z\right)\right|\\
	&\le C\{n^{-1/8}(\log(2pn/\gamma))^{7/8}+\gamma\},
	\end{split}
	\end{equation}
	where $C>0$ is a positive constant depending on $c_1$, $C_1$ and $C_2$ only.
\end{theorem}

\begin{remark}
	Notice that $\gamma$ can take all values in $(0,1)$ and not larger than the first term in the right-hand side of \eqref{e:ULB}. 
	Taking $\gamma=n^{-1/4}$, the right-hand side of~\eqref{e:ULB} will go to $0$ when $n,p\rightarrow\infty$ obeying $p\le e^{o(n^{1/7})}$. 
\end{remark}

\begin{remark}
	Theorem~\ref{thm:GauAppr} has a weaker assumption that the fourth moment exists comparing with Theorem~\ref{thm:TailProbMDT}, which requires a finite exponential moment. 
	Moreover, Theorem~\ref{thm:GauAppr} obtains a uniformly bound for the Gaussian approximation.
	But these relaxed conditions pay a cost on the convergence rate. 
	Although Theorem \ref{thm:TailProbMDT} requires a uniformly finite exponential moment, and it is only suitable for all $z\in [C\sqrt{\log{p}}, o(n^{1/6}(\log{p})^{-1/3}))$, it gives a more delicate estimation. 
\end{remark}


In the end, by using Theorem~\ref{thm:GauAppr} with $\gamma=n^{-1/4}$ we conclude in Corollary~\ref{coro:SM} that $\penaltyest_2$ is also a good approximation of $\lambda$. 
\begin{corollary}\label{coro:SM}
	Suppose that $\nabla L(\target|\designmatrix,\outcome)$ can be written as the following form
	\begin{equation*}
		\nabla L(\target|\designmatrix,\outcome)=\frac1n\sum_{i=1}^n\varWi,
	\end{equation*}
	with $\varWi=\fcnf(\designrowi,\outcomesi)\in\R^p$ being independent random vectors. 
	Suppose $c_1\le M_2^2\le C_1$ and $\E_n(\E(\max_{j\in [p]} \varWij^4))\le C_2$ with some positive constants $c_1$, $C_1$ and $C_2$.
	We have 
	\begin{equation*}
		\PP\left(c\|\nabla L(\target|\designmatrix,\outcome)\|_{\infty}\le\penaltyest_2\right)\ge 1-\alpha-O(n^{-1/8}(\log p)^{7/8}), 
	\end{equation*}
	where $\penaltyest_2=c(\sqrt{n})^{-1}z_{1-\alpha}$ with $c>1$, $\alpha\in(0,1)$ and $z_{1-\alpha}$ satisfying 
	\begin{equation*}
		\PP\left(\max_{j\in[p]\}}\left|\frac{1}{\sqrt{n}}\sum_{i=1}^{n}\varZij\right|\le z_{1-\alpha}\right)=1-\alpha. 
	\end{equation*} 
	If further assume $p\le e^{o(n^{1/7})}$, we have 
	\begin{equation*}
		\PP\left(c\|\nabla L(\target|\designmatrix,\outcome)\|_{\infty}\le\penaltyest_2\right)\ge 1-\alpha(1+o(1)).
	\end{equation*}
\end{corollary}

\begin{rem}
	Be similar to Remark~\ref{rem:coroMDT}, $\varWi$ in Corollary~\ref{coro:SM} can be written as $\varWi=\designrowi\epsilon_i$, if the function $L$ is the negative log-likelihood of generalized linear model, and suppose $\outcomesi=\mu(\designrowi'\target) + \epsilon_i$ for all $i=1,\ldots,n$ with $\mu(\designrowi'\target)=g^{-1}(\designrowi'\target)$. 
	Only if we assume $\epsilon_i$ satisfies $\E\epsilon_i=0$,  $\E\epsilon_i^2=\sigma^2<\infty$, and $\E_n(\max_{j\in [p]} \designsij^4\E\epsilon_i^4)\le C_2$ with some positive constant $C_2$, all the results of Corollary~\ref{coro:SM} hold. 
\end{rem}

\section{Examples and simulations}\label{sec:Applications}
In this section, we utilize the two methods in the previous section to estimate penalty levels of three types of $\ell_{1}$ penalized regressions including lasso \citep{Tibshirani1996Regression}, square-root lasso \citep{Belloni2011Square-root} and weighted $\ell_1$ penalized Poisson regression \citep{Jia2019Sparse}.

Recall the generalized linear model in Section\ref{sec:Introduction}. 
Suppose $(\designrow_{1},\outcomes_{1}),\cdots,$ $(\designrow_{n},\outcomes_{n})$ are independent pairs of observed data which are realizations of random vectors $(\varX_{1},\varY_{1}),\cdots,(\varX_{n},\varY_{n})$,
with covariates~$\varX_i\in\R^p$ and univariate response variables $\varY_{i}\in \R$ for all $i\in [n]$. 
$(\varX_i,\varY_{i})$ are assumed to satisfy the conditional distribution
\begin{equation}\label{e:glm}
\varY_i|\varX_i=\designrowi\sim F{\rm\ \ \ \ with}\ \ \ \ g(\E(\varY_i|\varX_i=\designrowi))=\designrowi'\target,
\end{equation}
where $F$ is some distribution in exponential family, $g(\cdot)$ a link function and $\target\in\mathbb{R}^{p}$ an unknown parameter vector to be estimated.
Denoting $\designrowi=(\designs_{i1},\cdots,\designs_{ip})'$, without loss of generality, we assume
$$\frac1{n}\sum\limits_{i=1}^n{\designsij}=0\ \ \ \ \ \ \ \ {\rm and}\ \ \ \ \ \ \ \frac1{n}\sum\limits_{i=1}^n{\designsij^2 }=1,\ \ \ \ \ \ \ \ \ \ \ {\rm for\ all}\ j\in[p].$$
Consider the following $\ell_{1}$-minimization problem:
\Be\label{e:estimator}
\est = \arg\min_{\parameter\in\R^p}\left\{ L(\parameter|\designmatrix,\outcome)+\lambda \|\parameter\|_1 \right\},
\Ee
where $L(\parameter|\designmatrix,\outcome)$, determined by the distribution of $\outcome|\designmatrix$, is assumed to be a convex function with respect to $\parameter$, and $\lambda>0$ is the penalty level.

In the ultrahigh dimensional regressions, for a (small) given $\alpha \in (0,1)$, we need to find some suitable $\lambda$ such that
\Bes
\PP\left(H \le \lambda\right)\ \ge\ 1-\alpha,
\Ees
where $H=c\|\nabla L(\target|\designmatrix,\outcome)\|_{\infty}$. 
The best choice of $\lambda$ is the $(1-\alpha)$-quantile of $H$. 
Unfortunately, $H$ is often not known especially when the errors are non-Gaussian, so $\lambda$ has to be chosen by heavy procedures such as cross-validation (CV).  

In what follows, we will use Corollary~\ref{coro:MDT} and Corollary~\ref{coro:SM} to approximate $\lambda$ in three specific examples, but not limited to them.
For simplicity, we use $L(\parameter)$ instead of $L(\parameter|\designmatrix,\outcome)$.

\subsection{Example 1: Lasso with linear models}\label{subsec:lasso}
Let $F$ be a distribution with mean $0$ and variance $\sigma^2<\infty$ and $g(\cdot)$ be the identity function. 
Then the model \eqref{e:glm} can be written as 
\begin{equation}\label{e:linearmodel}
\outcomesi=\designrowi'\target+\sigma\epsilon_i {\rm\ \ \ \ with}\ \ \ \ \sigma\epsilon_i \sim F.
\end{equation}

\noindent For Lasso \citep{Tibshirani1996Regression}, $L(\parameter)$ has the form $L(\parameter) =\sum_{i=1}^n(\outcomesi-\designrowi'\parameter)^2/2n$. 
Then,
\begin{equation*}
	\nabla L(\target) = -\frac1n \sum_{i=1}^{n}\designrowi (\outcomesi -\designrowi'\target)= -\frac{1}{n} \sum_{i=1}^n \sigma \designrowi\epsilon_i.
\end{equation*}

We firstly apply Corollary~\ref{coro:MDT} with $\varWi = -\sigma \designrowi\epsilon_i$, and it gives an approximation of $\lambda$ with probability nearly~$1-\alpha$.  
Observe that $\E_n(\E \varWij^2)=\theta^2=\sigma^2<\infty$. 
We assume that $\sup_{i \in [n], j \in [p]} \E e^{t |\sigma \designsij \epsilon_{i}|}<\infty$ for some $t \in (0,\infty)$. 
Then, defining
$$\penaltyest^{\operatorname{L}}_1(1-\alpha) = c\sigma(\sqrt{n})^{-1}\Phi^{-1}(1-\frac{\alpha}{2p}),$$
by \eqref{ineq:lambda1} we have
\begin{equation*}
	\PP\left(c\|\nabla L(\target)\|_\infty\le \penaltyest_1^{\operatorname{L}}(1-\alpha)\right) \ge  1-\alpha(1+O((\log p)^{5/2}/\sqrt{n})).
\end{equation*}
So, $\penaltyest_1^{\operatorname{L}}(1-\alpha)$ is a good approximation of $\lambda$ when $n, p$ are large.

Secondly, we apply Corollary~\ref{coro:SM} with $\varWi =-\sigma \designrowi\epsilon_i$, and it gives us another approximation for $\lambda$ under the conditions $M^{2}_{2}=\max_{j\in [p]}\E_{n}\left(\E\varWij^2\right) =\sigma^2<\infty$ and $\E_{n} \left[\E\left(\sigma^4 \max_{j \in [p]} \designsij^4 \epsilon_i^4\right)\right]<C_{2}$ for all $i \in [n]$ with some constant $C_2>0$. 
Let $e_i$ be i.i.d. standard Gaussian random variables, and let $z_{1-\alpha}^{\operatorname{L}}$ satisfy that
\begin{equation*}
	\PP\left(\max_{j\in [p]} \frac{1}{\sqrt{n}} \left|\sum_{i=1}^n \sigma \designsij e_i\right|\le z_{1-\alpha}^{\operatorname{L}}\right) = 1-\alpha.
\end{equation*}
Denote
$$\penaltyest^{\operatorname{L}}_2(1-\alpha) = c (\sqrt{n})^{-1}z_{1-\alpha}^{\operatorname{L}}.$$
Then, by Corollary~\ref{coro:SM}, we have
\begin{equation*}
	\PP\left(c\|\nabla L(\target)\|_\infty\le \penaltyest_2^{\operatorname{L}}(1-\alpha)\right) \ge  1-\alpha-O(n^{-1/8}(\log p)^{7/8}).
\end{equation*}

In a word, the analyses above show that given a confidence level $1-\alpha$, the two  choices $\penaltyest^{\operatorname{L}}_{1}(1-\alpha)$ and $\penaltyest^{\operatorname{L}}_2(1-\alpha)$ are good approximations of $\lambda$ for lasso. 
Additionally, they achieve high computational efficiency, which is shown in the following simulations. 
We also show the prediction errors under these two penalty levels. 

\begin{figure}
	\centering
	\includegraphics[scale = 0.35]{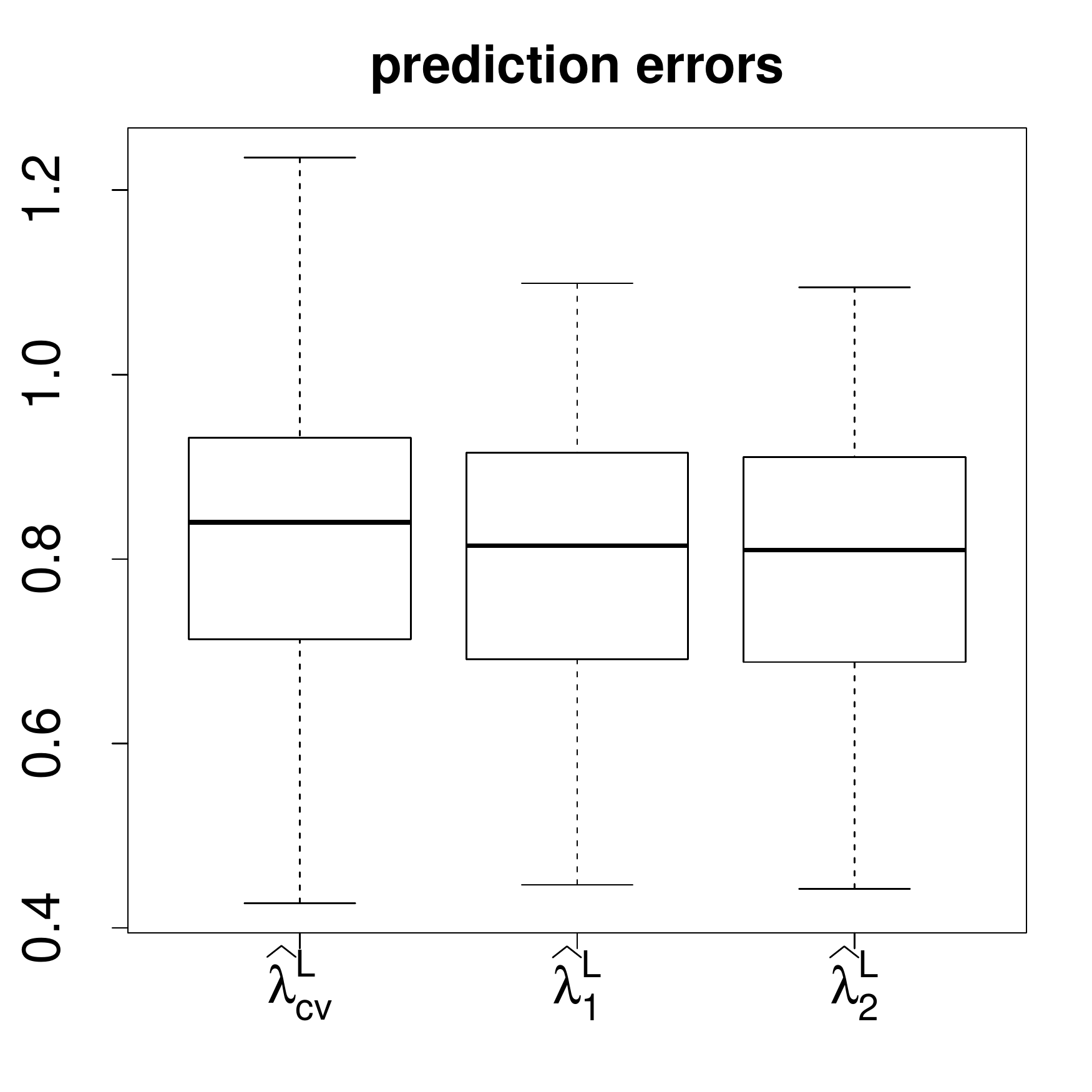}
	\caption{Prediction errors via three penalty selection methods for lasso with $n = 200$, $p = 1000$. 
		$\penaltyest_{\operatorname{cv}}^{\operatorname{L}}$: selected by the 10-fold CV; $\penaltyest_1^{\operatorname{L}}$: estimated based on the moderate deviation theorem; $\penaltyest_2^{\operatorname{L}}$: estimated based on Stein's method. }\label{fig:lasso}
\end{figure}

Under the model \eqref{e:linearmodel}, we use the following data settings: $n = 200$, $p = 1000$ and $\sigma = 1$. 
Let $\designrowi$ be generated from a $p$-dimensional normal distribution~$N(\boldsymbol{0}_p,\boldsymbol{\Sigma})$ with $\Sigma_{ij}=0.5^{|i-j|}$. 
The true parameter vector $\target$ just has 10 nonzero components. 
Let the first 10 components of $\target$ be non-zero and each non-zero component takes value randomly from $[-1,1]$. 
In practice, $\sigma$ is usually unknown, but there many ways to estimate it, such as the scaled lasso~\citep{Sun2012Scaled}. 
For convenience, we assume that $\sigma$ is known. 
We show the running times to obtain the $\penaltyest_1^{\operatorname{L}}(1-\alpha)$ and $\penaltyest_2^{\operatorname{L}}(1-\alpha)$ with $\alpha = 0.1, c = 1.01$ by simulations. 
For comparison, we also present the running time of obtaining penalty level by the 10-fold CV. 
We repeat each simulation 100 times and get the totally running time. 
The results are shown in Table~\ref{tab:lasso}. 
At the same time, we compare the prediction errors, defined as 
\begin{equation*}
	\operatorname{prediction\  error}=\sqrt{\frac{1}{n}\sum_{i=1}^n(\designrowi'(\est-\target))^2},
\end{equation*}
under these three ways. 
The results are shown in Figure~\ref{fig:lasso}.

From Table~\ref{tab:lasso}, we see that our two estimation methods require less than 1 min for 100 simulations, but the 10-fold CV needs more than 7 mins. 
Figure~\ref{fig:lasso} shows that the prediction errors in our two ways outperform that in the 10-fold CV.

\begin{table}
	\centering
	\caption{The running time to obtain each approximated penalty levels in 100 times repeations for lasso. (unit: second)}\label{tab:lasso}
	\begin{tabular}{ccccccccc}
		\hline 
		$c$&& $\alpha$ && 10-fold CV && $\penaltyest_1^{\operatorname{L}}(1-\alpha)$ && $\penaltyest_2^{\operatorname{L}}(1-\alpha)$\\
		1.01 && 0.1 && 450.700 && 0.005 && 57.185 \\
		\hline
	\end{tabular}
\end{table}

\subsection{Example 2: Square-root lasso under linear models}\label{subsec:sr-lasso}
We also use the linear model \eqref{e:linearmodel} in Section~\ref{subsec:lasso}.
For square-root Lasso \citep{Belloni2011Square-root}, $L(\parameter)$ has the form $L(\parameter) =\sqrt{\sum_{i=1}^n(\outcomesi-\designrowi'\target)^2/n}$. 
Then,
\begin{equation*}
	\nabla L(\target) = -\frac1n \sum_{i=1}^{n}\frac{\designrowi (\outcomesi -\designrowi'\target)}{\sqrt{\sum_{i=1}^n(\outcomesi-\designrowi'\target)^2/n}}= -\frac{1}{n} \sum_{i=1}^n \frac{\designrowi\epsilon_i}{\sqrt{\sum_{i=1}^n \epsilon_i^2/n}}.
\end{equation*}

Firstly, we use Corollary~\ref{coro:MDT} with $\varWi = -{\designrowi\epsilon_i}/{\sqrt{\sum_{i=1}^n \epsilon_i^2/n}}$ to estimate $\lambda$. 
Observe that $\E_n(\E \varWij^2)=\theta^2=1<\infty$. 
Given $\alpha\in(0,1)$, define
\begin{equation*}\label{e:sr-lamb1}
	\penaltyest^{\operatorname{sr}}_1 (1-\alpha)= c (\sqrt{n})^{-1}\Phi^{-1}(1-\frac{\alpha}{2p}).
\end{equation*}
Assume $\sup_{i \in [n], j \in [p]} \E e^{t|\designsij\epsilon_i|/\sqrt{\sum_{i=1}^n \epsilon_i^2/n}}<\infty$ for some $t \in (0,\infty)$. 
By Corollary~\ref{coro:MDT}, we have
\begin{equation*}\label{ineq:sr-1}
	\PP\left(c\|\nabla L(\target)\|_\infty\le \penaltyest_1^{\operatorname{sr}}(1-\alpha)\right)\le 1-\alpha(1+O((\log p)^{5/2}/\sqrt{n})).
\end{equation*}
Then, $\penaltyest_1^{\operatorname{sr}}(1-\alpha)$ is a good approximation of $\lambda$ when $n, p$ are large. 

Secondly, we apply Corollary~\ref{coro:SM} with $\varWi =-{\designrowi\epsilon_i}/{\sqrt{\sum_{i=1}^n \epsilon_i^2/n}}$ to find another approximation of $\lambda$. 
Let $e_i$ be i.i.d. standard Gaussian random variables and $z_{1-\alpha}^{\operatorname{sr}}$ satisfy that
\begin{equation*}\label{e:zalphasr}
	\PP\left(\max_{j\in [p]} \frac{1}{\sqrt{n}} \left|\sum_{i=1}^n \frac{\designsij e_i}{\sqrt{\sum_{i=1}^n e_i^2/n}}\right|\le z_{1-\alpha}^{\operatorname{sr}}\right) = 1-\alpha.
\end{equation*}
Denote
\begin{equation*}\label{e:sr-lamb2}
	\penaltyest^{\operatorname{sr}}_{2}(1-\alpha) = c(\sqrt{n})^{-1}z_{1-\alpha}^{\operatorname{sr}}. 
\end{equation*}
Under the conditions $M^{2}_{2}=\max_{j\in [p]}\E_{n}\left(\E\varWij^2\right) =1<\infty$ and $\E_{n} \left[\E\left(\max_{j \in [p]} \designsij^4 \epsilon_i^4/(\sum_{i=1}^n \epsilon_i^2/n)^2 \right)\right]<C_{2}$ for all $i \in [n]$ with some constant $C_2>0$, Corollary~\ref{coro:SM} implies that 
\begin{equation*}\label{ineq:GA-L}
	\PP\left(c\|\nabla L(\target)\|_\infty\le \penaltyest_2^{\operatorname{sr}}(1-\alpha)\right) \ge  1-\alpha-O(n^{-1/8}(\log p)^{7/8}).
\end{equation*}

Hence, $\penaltyest^{\operatorname{sr}}_{1}(1-\alpha)$ and $\penaltyest^{\operatorname{sr}}_2(1-\alpha)$ are two good approximations of $\lambda$ for square-root lasso under confidence level $1-\alpha$. 

\begin{figure}
	\centering
	\includegraphics[scale = 0.35]{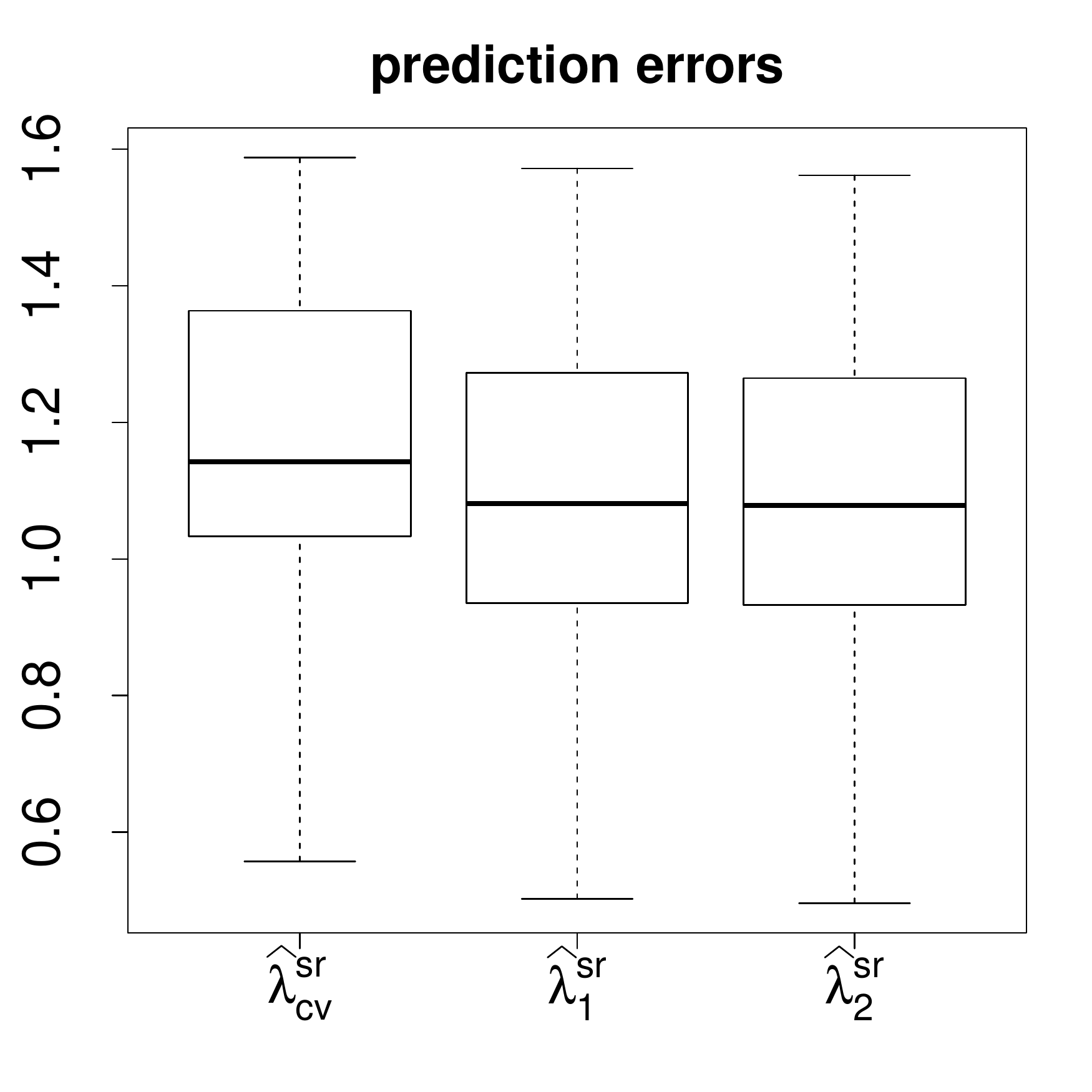}
	\caption{Prediction errors via three penalty selection methods for square-root lasso with $n = 200$, $p = 1000$. 
		$\penaltyest_{\operatorname{cv}}^{\operatorname{sr}}$: selected by the 10-fold CV; $\penaltyest_1^{\operatorname{sr}}$: estimated based on the moderate deviation theorem; $\penaltyest_2^{\operatorname{sr}}$: estimated based on Stein's method. }\label{fig:sr}
\end{figure}

Being similar to Example 1, we do some simulations to show the running time and prediction errors of square-root lasso under the three ways of selecting $\lambda$—$\penaltyest^{\operatorname{sr}}_{1}$,  $\penaltyest^{\operatorname{sr}}_2$, and 10-fold CV. 
The model and data settings are same as Example 1. 
The results are shown in Table~\ref{tab:sr} and Figure~\ref{fig:sr}. 
We see in Table~\ref{tab:sr} that the 10-fold CV costs more than 80 times time than our two ways. 
More than that, the prediction errors under our two selections outperforms that under the 10-fold CV. 

\begin{table}
	\centering
	\caption{The running time to obtain each approximated penalty levels in 100 times repeations for square-root lasso. (unit: second)}\label{tab:sr}
	\begin{tabular}{ccccccccc}
		\hline
		$c$&& $\alpha$ && 10-fold CV && $\penaltyest_1^{\operatorname{sr}}(1-\alpha)$ && $\penaltyest_2^{\operatorname{sr}}(1-\alpha)$\\
		\\
		1.01 && 0.1 && 4188.876  && 0.001 &&  48.597\\
		\hline
	\end{tabular}
\end{table}

\subsection{Example 3: Poisson regression with $\ell_1$ penalized weighted score function method}
In this section, we apply the Corollary~\ref{coro:MDT} and Corollary~\ref{coro:SM} into Poisson regression case.  
Firstly we introduce the Poisson regression with $\ell_1$ penalized weighted score function method (LPWSF) \citep{Jia2019Sparse}. 
This method is different from the traditional $\ell_1$ penalized maximum log-likelihood estimation \citep{Li2015Consistency}. 
As the name implies, it added weight on the score function of the Poisson distribution.

Suppose $F$ is Poisson distribution with parameter $\mu(\designrowi)$ and let link function $g(x)=\log{x}$. 
Then \eqref{e:glm} becomes
\begin{equation}\label{e:modelPoi}
\varY_i|\varX_i=\designrowi\sim {\rm Poisson}(\mu(\designrowi)){\rm\ \ \ \ with}\ \ \ \ \log(\mu(\designrowi))=\designrowi'\target,
\end{equation}
where $\target\in\mathbb{R}^{p}$ is an unknown parameter vector to be estimated.
The Poisson estimator obtained by LPWSF method is defined by \eqref{e:estimator} with
$$L(\parameter) = \frac1{n}\sum\limits_{i=1}^{n}2(\outcomesi e^{-\designrowi'\parameter/2}+e^{\designrowi'\parameter/2}).$$
Then, the gradient of $L(\parameter)$ valued at $\parameter = \target$ has the form
$$\nabla L(\target)= -\frac1{n}\sum\limits_{i=1}^{n}\frac{\designrowi(\outcomesi -e^{\designrowi'\target})}{\sqrt{e^{\designrowi'\target}}}.$$
Let $\epsilon_i={\outcomesi -e^{\designrowi'\target}}/{\sqrt{e^{\designrowi'\target}}}$ and then wer observe that $\E \epsilon_i =0$ and $\E \epsilon_i^2 = 1$. 

Firstly, we apply Corollary~\ref{coro:MDT} with $\varWi=-\designrowi\epsilon_i$  to find an approximation of $\lambda$. 
Observe $\E_n(\E \varWij^2)=\theta^2=1<\infty$ and assume $\sup_{i \in [n], j \in [p]} \E e^{t |\designsij \epsilon_{i}|}<\infty$ for some $t \in (0,\infty)$. 
Denote 
\begin{equation*}\label{e:poi-lamb1}
	\penaltyest^{\operatorname{P}}_{1}(1-\alpha) =c (\sqrt{n})^{-1}\Phi^{-1}(1-\frac{\alpha}{2p}). 
\end{equation*}
By \eqref{ineq:lambda1}, we have
\begin{equation*}
	\PP\left(c\|\nabla L(\target)\|_\infty\le \penaltyest^{\operatorname{P}}_{1}(1-\alpha)\right)\ge 1-\alpha(1+O((\log p)^{5/2}/\sqrt{n})).
\end{equation*}
Then, $\penaltyest_1^{\operatorname{sr}}(1-\alpha)$ is a good approximation of $\lambda$ when $n, p$ are large. 

Secondly, using a Gaussian approximation in Corollary~\ref{coro:SM}, we define 
\begin{equation*}\label{e:poi-lamb2}
	\penaltyest^{\operatorname{P}}_2(1-\alpha) = c(\sqrt{n})^{-1}z_{1-\alpha}^{\operatorname{P}},
\end{equation*}
where $z_{1-\alpha}^{\operatorname{P}}$ satisfies 
$$\PP\left(\max_{j\in [p]} \frac{1}{\sqrt{n}} \left|\sum_{i=1}^n \designsij e_i\right|\le z_{1-\alpha}^{\operatorname{P}}\right) = 1-\alpha,$$
with $e_i$ being independent standard normal random variables. 
Then, under the conditions $M^{2}_{2} =\max_{j\in [p]}\E_{n}\left(\E\varWij^2\right) =1$ and $\E_{n} \left[\E\left(\max_{j \in [p]} \designsij^4 \epsilon_i^4\right)\right]<C_{2}$ for all $i \in [n]$ with some constant $C_2>0$, Corollary~\ref{coro:SM} yields that 
\begin{equation*}
	\PP\left(c\|\nabla L(\target)\|_\infty\le \penaltyest^{\operatorname{P}}_{2}(1-\alpha)\right)\ge 1-\alpha-O(n^{-1/8}(\log p)^{7/8}).
\end{equation*}

Hence, for weighted $\ell_1$ penalized Poisson regression, $\penaltyest^{\operatorname{P}}_1(1-\alpha)$ and $\penaltyest^{\operatorname{P}}_2(1-\alpha)$ are two suitable choices of $\lambda$. 
We conduct some simulations to further confirm the results below.

\begin{figure}
	\centering
	\includegraphics[scale = 0.35]{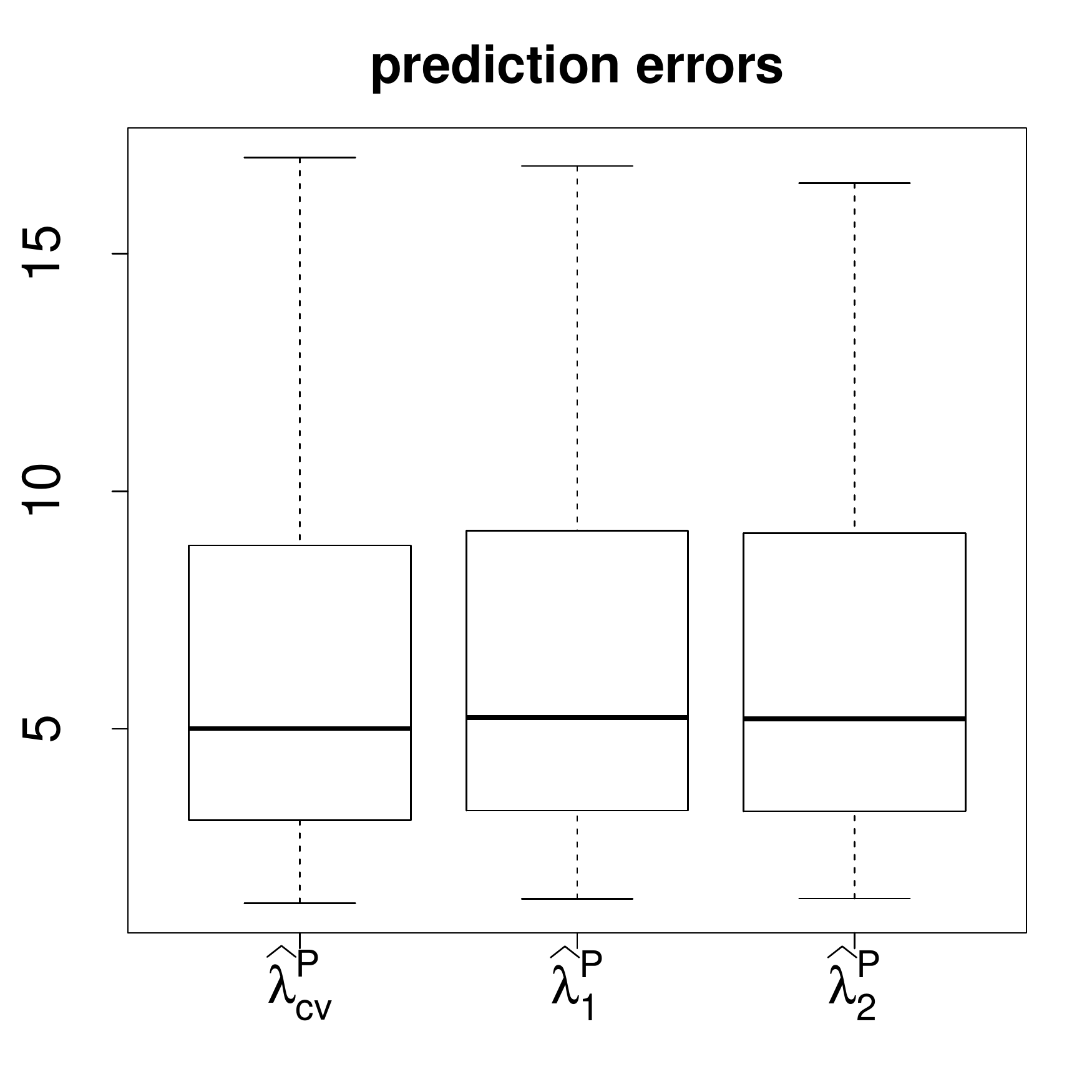}
	\caption{Prediction errors via three penalty selection methods for Poisson regression with LPWSF method ($n = 200$, $p = 1000$). 
		$\penaltyest_{\operatorname{cv}}^{\operatorname{P}}$: selected by the 10-fold CV; $\penaltyest_1^{\operatorname{P}}$: estimated based on the moderate deviation theorem; $\penaltyest_2^{\operatorname{P}}$: estimated based on Stein's method. }\label{fig:poi}
\end{figure}

Being similar to the previous two examples, we also compare the running times of three ways to select the penalty level and exhibit the prediction errors of each method via box-plot. 
Under the model \eqref{e:modelPoi}, all the settings of $n,p$, $\designrowi$, $\target$ and so on are same with lasso except that $\outcomesi $ is generated from a Poisson distribution with parameter $e^{\designrowi'\target}$ for each $i\in[n]$. 
We repeat the simulations 100 times for each method. 
The running time is shown in Table~\ref{tab:poi}, which
shows that the 10-fold CV costs more than 26 times time than our two ways to estimate $\lambda$. 
We get a similar prediction error behavior with the 10-fold CV, which are shown in Figures~\ref{fig:poi}. 
So, our two ways are more competitive than the 10-fold CV. 

\begin{table}[H]
	\centering
	\caption{The running time to obtain each approximated penalty level in 100 times repetitions for the Poisson regression with the LPWSF method. (unit: second)}\label{tab:poi}
	\begin{tabular}{ccccccccc}
		\hline
		$c$ && $\alpha$ && 10-fold CV && $\penaltyest_1^{\operatorname{P}}(1-\alpha)$ && $\penaltyest_2^{\operatorname{P}}(1-\alpha)$\\
		\\
		1.01  && 0.1 && 1121.308  && 0.003  && 41.762 \\
		\hline
	\end{tabular}
\end{table}

\section{Conclusion}\label{sec:conclusion}
We proposed two theoretical approximations for the penalty level of $\ell_{1}$ penalized generalized linear regressions. 
The main skills were truncation technique, moderate deviation theorem, and Stein's method.
These skills can be also used for any other applications, which need a Gaussian approximation. 
We applied our approximated penalty levels to three types of high-dimensional $\ell_{1}$ penalized regressions, lasso, square-root lasso and Poisson regression with $\ell_{1}$ penalized weighted score function method. 
The simulation results showed that our approximated penalty levels produced comparable accuracy on the prediction error comparing with 10-fold CV and also achieved high computational efficiency. 

In addition, the two theoretical approximations are not limited in the generalized linear models. 
They are suitable for all the models only if the gradient of their objective function $L$ can be written as a average of $n$ independent random vectors (see more details in Corollary~\ref{coro:MDT} and Corollary~\ref{coro:SM}).

\paragraph{Acknowledgements}
We would like to gratefully thank Prof. Jian-feng Yao for stimulating discussions at a conference in ShenZhen.

\appendix


\section{Proof of Theorem \ref{thm:TailProbMDT}}\label{appendix:MDT}
Denote $\hat{x}_{ij} = \varWij1_{\{|\varWij|\le A\}}$ and $\check{x}_{ij}= \varWij1_{\{|\varWij|> A\}},$
where $A$ will be chosen later.
Observing that $\varWij = \hat{x}_{ij}+\check{x}_{ij}-\E \hat{x}_{ij}-\E \check{x}_{ij}$ we have for $z\geq 0$, 
\begin{equation}\label{ineq:I}
\begin{split}
&\PP\left(\max_{1 \le j \le p} |\sumWj| > z\right)\\
&\le \sum_{j=1}^p \PP\left(\left|\frac{1}{\sqrt{n}}\sum_{i=1}^n \varWij\right| > z\right)\\
&= \sum_{j=1}^p \PP\left(\frac1{\sqrt{n}}\left|\sum_{i=1}^n \left(\hat{x}_{ij}+\check{x}_{ij}-\E \hat{x}_{ij}-\E \check{x}_{ij}\right)\right|>z\right)\\
&\le \sum_{j=1}^p \left\{\PP\left(\frac1{\sqrt{n}}\left|\sum_{i=1}^n\left(\hat{x}_{ij}-\E \hat{x}_{ij}-\E \check{x}_{ij}\right)\right|>z\right)+\PP\left(\sup_{i\in[n]}|\varWij|>A\right)\right\}\\
&\le \sum_{j=1}^p \left\{\PP\left(\left|\sum_{i=1}^n\left(\hat{x}_{ij}-\E \hat{x}_{ij}\right)\right|>\sqrt{n}z-\left|\sum_{i=1}^n \E \check{x}_{ij}\right|\right)+\PP\left(\sup_{i\in[n]}|\varWij|>A\right)\right\}\\
&= \sum_{j=1}^p \left(I_j^1 + I_j^2\right),
\end{split}
\end{equation}
where
\begin{align*}
	I_j^1 = \PP\left(\left|\sum_{i=1}^n\left(\hat{x}_{ij}-\E \hat{x}_{ij}\right)\right|>\sqrt{n}z-\left|\sum_{i=1}^n \E \check{x}_{ij}\right|\right) \ {\rm and }\  I_j^2 = \PP\left(\sup_{i\in[n]}|\varWij|>A\right).
\end{align*}

Firstly, we estimate $I_j^2$. For each $j\in[p]$, by exponential Chebyshev's inequality and condition $\sup_{i\in[n],j\in[p]}\E e^{t_1|\varWij|}<\infty$ for some $t_1>0$, we have
\begin{equation}\label{ineq:I2j}
I_j^2\le \sum_{i=1}^n \PP \left(|\varWij|>A\right)\le n e^{-t_1A}\E e^{t_1|\varWij|}\le Cne^{-t_1A},
\end{equation}
where $t_1>0$ is some constant.

Before estimating $I_j^2$, we make some preparations. Observe that
\begin{align*}
	|\mathbb{E} \check{x}_{ij}|&\le\mathbb{E}|\check{x}_{ij}|=\mathbb{E}|\varWij|1_{\{|\varWij|>A\}}=
	\int_{A}^{+\infty}z d F(z)+\int_{-\infty}^{-A}-z d F(z) \\
	&=\left\{z(F(z)-1)|_{A}^{+\infty}-\int_{A}^{+\infty}(F(z)-1)dz\right\}+\left%
	\{\int_{-\infty}^{-A}F(z)dz-zF(z)|_{-\infty}^{-A}\right\} \\
	&\le A(1-F(A))+\int_{A}^{+\infty}Ce^{-t_1z}dz +
	\int_{-\infty}^{-A}Ce^{t_1z}dz+AF(-A) \\
	&\le C(A+2/t_1)e^{-t_1A}.
\end{align*}
Denote $u=C(A+2/t_1)e^{-t_1A}$ and then $|\mathbb{E} \check{x}_{ij}|\le u$. Let $\eta_{ij} = {(\hat{x}_{ij}-\E \hat{x}_{ij})}/{2A}$. Then, $\E \eta_{ij}=0$ and $|\eta_{ij}|\le 1$ for all $i\in [n], j\in[p]$. Denote $\sigma_{nj}^2=\sum_{i=1}^n\E \eta_{ij}^2$ and $T_{nj}= \sum_{i=1}^n \E |\eta_{ij}|^3/\sigma_{nj}^3$. By calculating we have
\begin{equation*}
	\begin{split}
		\sigma_{nj}^2&=\sum_{i=1}^n\E\left(\frac{\hat{x}_{ij}-\E \hat{x}_{ij}}{2A}\right)^2\le \frac1{4A^2}\sum_{i=1}^n\E \varWij^2 = \frac{n\theta^2}{4A^2}\\
		T_{nj}&=\sum_{i=1}^n  \E |\eta_{ij}|^3/\sigma_{nj}^3\le \sum_{i=1}^n  \E |\eta_{ij}|^2/\sigma_{nj}^3 = \frac{1}{\sigma_{nj}}.
	\end{split}
\end{equation*}
For convenient, we denote $T_{nj} = O(1)A\left(\sqrt{n}\right)^{-1}$.
Then using Lemma \ref{lem:MDT} we have
\begingroup
\allowdisplaybreaks
\begin{align*}
	I_j^1&\le \PP\left(\left|\sum_{i=1}^n\left(\hat{x}_{ij}-\E \hat{x}_{ij}\right)\right|>\sqrt{n}z-\left(\sum_{i=1}^n \left|\E \check{x}_{ij}\right|^2\right)^{1/2}\right)\\
	&\le \PP\left(\left|\sum_{i=1}^n\left(\hat{x}_{ij}-\E \hat{x}_{ij}\right)\right|>\sqrt{n}(z-u)\right)\\
	&=\PP\left(\left|\sum_{i=1}^n\eta_{ij}\right|>\frac{\sqrt{n}(z-u)}{2A\sigma_{nj}}\sigma_{nj}\right)\\
	&\le\PP\left(\left|\sum_{i=1}^n\eta_{ij}\right|>(z-u)\sigma_{nj}\right)\\
	&= 2\left(1+O(1)\left(z-u\right)^3 T_{nj}\right)\bar{\Phi}\left(z-u\right)\\
	&\le2\left(1+O(1)\frac{A(z-u)^3}{\sqrt{n}}\right)\left(\bar{\Phi}\left(z\right)+\frac{u}{\sqrt{2\pi}}\right)\\
	&= 2\bar{\Phi}\left(z\right)\left(1+O(1)\frac{A(z-u)^3}{\sqrt{n}}\right)\left(1+\frac{u}{\sqrt{2\pi}\bar{\Phi}\left(z\right)}\right).
\end{align*}
\endgroup
Combining \eqref{ineq:I}, \eqref{ineq:I2j} and the inequality above, we have
\begin{align*}
	&\PP\left(\max_{1 \le j \le p} |\sumWj| \le z\right) \\
	&= 1-\PP\left(\max_{1 \le j \le p} |\sumWj| > z\right)\\
	&\ge 1-2p\bar{\Phi}\left(z\right)\left(1+O(1)\frac{A(z-u)^3}{\sqrt{n}}\right)\left(1+\frac{u}{\sqrt{2\pi}\bar{\Phi}\left(z\right)}\right) - Cnpe^{-t_1 A}.
\end{align*}
Taking $A = {3\log{p}}/{t_1}$ and $u = {O(1)\log{p}}/{p^3}$, we obtain the desired result. \qed

\section{Proof of Theorem \ref{thm:GauAppr}}\label{appendix:SM}
The following two lemmas, which are Theorem 2.2 and Lemma 2.2 of \citep{chernozhukov2013Gaussian} respectively, are the keys to obtain this theorem.
\begin{lem}\label{lem:chernozhukov2013}
	Suppose that there are some constants $0<d_1<D_1$ such that $d_1\le \E_n( \E(\varWij^2))\le D_1$ for all $j\in[p]$. Then for every $\gamma\in (0,1)$,
	\Bes
	\begin{split}
		&\quad\sup_{z\in \mathbb{R}}\left|\PP\left(\max_{1 \le j \le p} \sumWj \le z\right)-\PP\left(\max_{1 \le j \le p} \sumZj \le z\right)\right|\\
		&\le D\{n^{-1/8}(M_3^{3/4}\vee M_4^{1/2})(\log(pn/\gamma))^{7/8}+n^{-1/2}(\log{pn/\gamma})^{3/2}a(\gamma)+\gamma\},
	\end{split}
	\Ees
	where $D>0$ is a constant that depends on $d_1$ and $D_1$ only.
\end{lem}

\begin{lemma}\label{lem:chernozhukov2013-2}
	Let $f:[0, \infty)\rightarrow [0, \infty)$ be a Young-Orlicz modulus, and let $f^{-1}$ be the inverse function of $f$. Let $B_1>0$ and $B_2>0$ be constants such that $(\E(\varWij^2))^{1/2}\le B_1$ for all $i\in[n], j\in[p]$ and $\E_n(\E[f(\max_{j\in [p]} |\varWij|/B_2)])\le 1$. Then under the condition of Lemma \ref{lem:chernozhukov2013},
	\begin{equation*}
		a(\gamma)\le D'\max\{B_2 f^{-1}(n/\gamma), B_1 \sqrt{\log(pn/\gamma)}\},
	\end{equation*}
	where $D'>0$ is a constant that depends on $d_1$ and $D_1$ only.
\end{lemma}

For all $i\in[n]$, define two $2p$-dimensional vectors $\varWitilde$ and $\varZitilde$ with $\varWitilde=(\varWi',-\varWi')'$ and $\varZitilde=(\varZi',-\varZi')'$. Then, $\sumWtilde$ and $\sumZtilde$ become two $2p$-dimensional vectors such that
\begin{align*}
	&\sumWtilde = \frac{1}{\sqrt{n}}\sum_{i=1}^n\varWitilde\ \ \ \ \ {\rm and}\ \ \ \ \ \sumWltilde = \frac{1}{\sqrt{n}}\sum_{i=1}^n\varWstilde_{il},\\
	&\sumZtilde = \frac{1}{\sqrt{n}}\sum_{i=1}^n\varZitilde\ \ \ \ \ {\rm and}\ \ \ \ \ \sumZltilde = \frac{1}{\sqrt{n}}\sum_{i=1}^n\varZstilde_{il},
\end{align*}
where $l\in[2p]$. Obviously, we have the following relations
\begin{equation}\label{e:2p-dim}
\max\limits_{1\le j \le p}|\sumWj|=\max\limits_{1\le l \le 2p} \sumWltilde,\ \ \ \ \ \ \ \ \max\limits_{1\le j \le p}|\sumZj|=\max\limits_{1\le l \le 2p} \sumZltilde.
\end{equation}
Since for all $j\in[p]$, $\varWstilde_{ij}=-\varWstilde_{i(j+p)} = \varWij$ and $\varZstilde_{ij}=-\varZstilde_{i(j+p)} = \varZij$, we have
\begin{equation}\label{e:moments}
\max_{l\in [2p]}\left(\E_{n}\left(\E\left|\varWstilde_{il}\right|^k\right)\right)^{1/k} = \max_{j\in [p]}\left(\E_{n}\left(\E\left|\varWij\right|^k\right)\right)^{1/k}=M_k.
\end{equation}

Let Young-Orlicz modulus $f$ in Lemma \ref{lem:chernozhukov2013-2} take power function $f(u) = u^{4}$ with inverse $f^{-1}(t) = t^{1/4}$. Since $\E_n(\E(\max_{j\in [p]} \varWij^4))\le C_2$, by Lemma \ref{lem:chernozhukov2013-2} we have
\begin{equation}\label{ineq:a(gamma)}
a(\gamma) \lesssim (n/\gamma)^{1/4}.
\end{equation}
Further by the relations \eqref{e:2p-dim}, \eqref{e:moments} and \eqref{ineq:a(gamma)}, we can get our desired result \eqref{e:ULB} from Lemma \ref{lem:chernozhukov2013} and Lemma \ref{lem:chernozhukov2013-2} straightforwardly. \qed



\bibliographystyle{elsarticle-harv} 
\bibliography{EstimatingPenaltyLevel}





\end{document}